\documentclass{article}

\usepackage{amsfonts, amsmath, amssymb, amsthm, amscd}
\usepackage{ascmac}

\usepackage{verbatim}

\usepackage{graphics}
\usepackage[all,2cell]{xypic}

\usepackage{mathptmx}

\usepackage{helvet}
\usepackage{courier}        
\usepackage{type1cm} 

\objectmargin+{1mm}
\labelmargin+{0.8mm}
\SelectTips{cm}{12}

\theoremstyle{plain}  
\newtheorem{thm}{Theorem}[section]

\newtheorem{cor}[thm]{Corollary}
\newtheorem{lem}[thm]{Lemma}
\newtheorem{prop}[thm]{Proposition}

\theoremstyle{definition}

\newtheorem{df}[thm]{Definition}

\newtheorem{para}[thm]{}

\theoremstyle{remark}

\DeclareMathOperator{\cA}{\mathcal{A}}
\DeclareMathOperator{\cB}{\mathcal{B}}
\DeclareMathOperator{\cC}{\mathcal{C}}
\DeclareMathOperator{\calD}{\mathcal{D}}

\DeclareMathOperator{\cF}{\mathcal{F}}

\DeclareMathOperator{\cS}{\mathcal{S}}

\DeclareMathOperator{\bA}{\mathbf{A}}
\DeclareMathOperator{\bB}{\mathbf{B}}
\DeclareMathOperator{\bC}{\mathbf{C}}

\DeclareMathOperator{\bbK}{\mathbb{K}}

\DeclareMathOperator{\bbN}{\mathbb{N}}

\DeclareMathOperator{\fF}{\mathfrak{F}}

\DeclareMathOperator{\fj}{\mathfrak{j}}

\DeclareMathOperator{\fS}{\mathfrak{S}}

\UseAllTwocells

\def\sn{\smallskip\noindent}
\def\mn{\medskip\noindent}

\def\enumidef{\renewcommand{\labelenumi}{$\mathrm{(\roman{enumi})}$}}
\def\enumiidef{\renewcommand{\labelenumii}{$\mathrm{(\roman{enumii})}$}}

\newcommand{\cf}{\textrm{cf.}\;}

\newcommand{\Cone}{\operatorname{Cone}}

\newcommand{\id}{\operatorname{id}}
\newcommand{\im}{\operatorname{Im}}

\newcommand{\isoto}{\overset{\scriptstyle{\sim}}{\to}}

\newcommand{\linf}{\leftarrowtail}

\newcommand{\Ob}{\operatorname{Ob}}
\newcommand{\onto}[1]{\stackrel{#1}{\to}}

\newcommand{\rdef}{\twoheadrightarrow}
\newcommand{\rinc}{\hookrightarrow}
\newcommand{\rinf}{\rightarrowtail}

\newcommand{\st}{\operatorname{st}}
\newcommand{\sst}{\operatorname{sst}}
\newcommand{\sus}{\operatorname{sus}}

\title{Delooping of the $K$-theory of 
strictly derivable Waldhausen categories}
\date{}

\author{Satoshi Mochizuki}

\begin{document}

\maketitle

\section*{Introduction}

In this short note, 
for a morphism of Waldhausen categories 
$f\colon\bA=(\cA,w_{\bA})\to \bB=(\cB,w_{\bB})$, 
we will define $\Cone f$ to be a Waldhausen category. 
There exists the canonical morphism of 
Waldhausen categories $\kappa_f\colon \bB\to \Cone f$ 
(see Definition~\ref{df:Cone of morphisms of Waldhausen categories}). 
We will show that the sequence $\bA\onto{f}\bB\onto{\kappa_f}\Cone f$ induces 
fibration sequence of spaces 
$K(\bA)\onto{K(f)}K(\bB)\onto{K(\kappa_f)}K(\Cone f)$ 
on connective $K$-theory 
(see Theorem~\ref{thm:Puppe sequences for connective K-theory}). 
Moreover we will define a notion of 
{\it strictly derivable Waldhausen categories} 
(see Definition~\ref{df:strictly derivable}) 
and define non-connective $K$-theory for strictly derivable Waldhausen categories (see Definition~\ref{df:nonconnective K-theory}). 

\section{Admissible classes of morphisms and $w$-thick subcategories}
\label{sec:adm class and Serre cat}

In this section, let 
$\cC$ be a category with cofibrations. 
We fix a zero object $0$ of $\cC$. 

\begin{df}[\bf Admissible class of morphisms]
\label{df:admissible class of morphisms}
We say that a class $w$ of morphisms in $\cC$ is 
{\it admissible} if 
$w$ satisfies the following three conditions. 
\begin{itemize}
\item[\ ]
{\bf (Closed under isomorphisms).}
$w$ contains all isomorphisms.
\item[\ ]
{\bf (Two out of three property for compositions).} 
$w$ satisfies the two out of three property. 
Namely for a pair of composable morphisms $x\onto{f}y\onto{g}z$, 
if two of $gf$, $f$ and $g$ are in $w$, then the third one is also in $w$. 

\item[\ ]
{\bf (Two out of three property for cofibration sequences).} 
In the commutative diagram of cofibration sequences below
\begin{equation}
\label{eq:cofib seq}
\xymatrix{
x \ar@{>->}[r]^j \ar[d]_a & y \ar@{->>}[r]^q \ar[d]_b & y/x \ar[d]^c\\
x' \ar@{>->}[r]_{j'} & y' \ar@{->>}[r]_{q'}  & y'/x', 
}
\end{equation}
if two of $a$, $b$ and $c$ are in $w$, 
then the third one is also in $w$.
\end{itemize}
\end{df}

\begin{lem}
\label{lem:admissible class satisfies the gluing axiom}
Assume that $\cC$ is enriched over the category of abelian groups. 
If $w$ is an admissible class of morphisms in $\cC$, then 
$w$ satisfies the gluing axiom of {\rm\cite{Wal85}}. 
In particular the pair $(\cC,w)$ is a category with 
cofibrations and weak equivalences.
\end{lem}

\begin{proof}
Since $\cC$ is close under finite coproducts, 
$\cC$ is an additive category. 
In the commutative diagram below, 
assume that $a$, $b$ and $c$ are in $w$. 
$$
{\footnotesize{
\xymatrix{
y\ar[d]_b & x \ar[l] \ar@{>->}[r] \ar[d]_a & z \ar[d]^c\\
y' & x' \ar[l] \ar@{>->}[r] & z'.
}}}
$$
Then there are commutative diagrams of cofibrations sequences.
$$
{\footnotesize{
\xymatrix{
y\ar@{>->}[r] \ar[d]_b & y\oplus z \ar@{->>}[r] \ar[d]_{{\tiny{\begin{pmatrix}b&0\\ 0 &c\end{pmatrix}}}} & z \ar[d]^c\\
y'\ar@{>->}[r] & y'\oplus z' \ar@{->>}[r] & z',
}\ \ \ \ \ 
\xymatrix{
x \ar@{>->}[r] \ar[d]_a & y\oplus z \ar[d]_{{\tiny{\begin{pmatrix}b&0\\ 0 &c\end{pmatrix}}}} \ar@{->>}[r] & y\sqcup_x z \ar[d]^{b\sqcup_a c}\\
x' \ar@{>->}[r] & y'\oplus z' \ar@{->>}[r] & y'\sqcup_{x'} z'.
}
}}
$$
By the left diagram above, it turns out that ${\tiny{\begin{pmatrix}b&0\\ 0 &c\end{pmatrix}}}$ is in $w$ and 
therefore $b\sqcup_a c$ is also in $w$ by the right diagram above. 
\end{proof}

We fix $w$ an admissible class of morphisms in $\cC$. 

\begin{df}[\bf Thick subcategory, $w$-closed subcategory]
\label{df:thick subcategory}
Let $\cS$ be a class of objects in $\cC$. 
We sometimes regard $\cS$ as a full subcategory of $\cC$. 
We say that $\cS$ is 
a {\it left prethick subcategory of $\cC$} 
if it contains the zero object $0$ of $\cC$ and if 
for a cofibration sequence $x\rinf y \rdef y/x$ in $\cC$ with 
$y/x\in\cS$, $x$ is in $\cS$ if and only if $y$ is in $\cS$. 
We say that $\cS$ is 
a {\it right prethick subcategory of $\cC$} 
if it contains the zero object $0$ of $\cC$ and if 
for a cofibration sequence $x\rinf y \rdef y/x$ in $\cC$ with 
$x\in\cS$, $y$ is in $\cS$ if and only if $y/x$ is in $\cS$. 
We say that $\cS$ is a {\it prethick subcategory of $\cC$} 
if it is both a left and right prethick subcategory of $\cC$.

We say that $\cS$ is {\it left thick} (respectively {\it right thick}, {\it thick}) if $\cS$ is left prethick (respectively right prethick, prethick) and moreover if $\cS$ is closed under retractions. 
The last condition means that for an object $x$ in $\cC$ if there exists 
morphisms $x\onto{i}y\onto{p}x$ such that $pi=\id_x$ and $y$ is in $\cS$, then 
$x$ is also in $\cS$. 
We say that $\cS$ is {\it $w$-closed} 
if it contains the zero object $0$ and if 
for an object $x$ in $\cC$, if there exists an object $y$ 
in $\cS$ and if there exists a zig-zag sequence of morphisms in $w$ 
which connects $x$ and $y$, then $x$ is also in $\cS$. 

For an admissible class $u$ of morphisms in $\cC$, 
we write $\cC^u$ for the full subcategory 
of $\cC$ consisting of those objects $x$ such that 
the canonical morphism $0\to x$ is in $u$. 
Then by the lemma below, $\cC^u$ is a thick subcategory 
of $\cC$. 

We say that $\cS$ is {\it left $w$-prethick} (respectively {\it right $w$-prethick}, {\it right $w$-thick}, {\it $w$-prethick}, 
{\it $w$-thick}) if it is left prethick 
(respectively left thick, right prethick, right thick, prethick, thick) and 
if it contains $\cC^w$. 
\end{df}

\begin{lem}
\label{lem:ex of w-Serre subcategory}
For an admissible class $u$ of morphisms in $\cC$, 
$\cC^u$ is a thick subcategory of $\cC$. 
Moreover if $u$ contains $w$, then 
$\cC^u$ is $w$-closed.
\end{lem}

\begin{proof}
For a cofibration sequence $x\rinf y \rdef y/x$, 
by applying the axiom of admissible classes to 
the commutative diagram below, 
it turns out that if two of $x$, $y$ and $y/x$ are in $\cC^u$, 
then the third one is also in $\cC^u$.
$$
{\footnotesize{
\xymatrix{
0 \ar@{>->}[r] \ar[d] & 0 \ar[d] \ar@{->>}[r] & 0 \ar[d]\\
x \ar@{>->}[r] & y \ar@{->>}[r] & y/x.
}}}
$$

Next let $x\onto{f}y$ be a morphism in $w$ and 
assume that $x$ (resp. $y$) is in $\cC^u$, then 
$y$ (resp. $x$) is also in $\cC^u$ 
by the two out of three property of $u$. 
$$
{\footnotesize{\xymatrix{
x \ar[rr]^f & & y\\
& 0\ar[lu] \ar[ru]. 
}}}
$$
\end{proof}

The rest of this section, 
we assume that the pair $(\cC,w)$ satisfies 
the following {\it factorization axiom} in \cite[A.5]{Sch06}.

\begin{para}[\bf Factorization axiom]
\label{df:factorization axiom}
For any morphism $x\onto{f}y$ in $\cC$, 
there exists a cofibration $i_f\colon x\rinf u_f$ and 
a morphism $q_f\colon u_f\to y$ in $w$ such that $f=q_fi_f$. 
In this case, 
we call a triple $(i_f,q_f,u_f)$ a ($w$-){\it factorization of $f$}. 
\end{para}

Recall that in \cite[1.1]{Cis10}, 
we say that a Waldhausen category is 
{\it derivable} if it satisfies the two out of three property for 
compositions and if it satisfies the factorization axiom. 
Thus we define the following notion.

\begin{df}[\bf Strictly derivable Waldhausen category]
\label{df:strictly derivable}
We say that a Waldhausen category $(\cC,w)$ is 
{\it strictly derivable} if it satisfies the 
factorization axiom and $w$ is an admissible class. 
\end{df}

The following two lemmata are useful to construct the correspondence 
between admissible classes and prethick subcategories.

\begin{lem}
\label{lem:func factorization}
Let 
\begin{equation}
\label{eq:commutative square}
\xymatrix{
x \ar[r]^i \ar[d]_a & y \ar[d]^b\\
x' \ar[r]_{i'} & y'
}
\end{equation}
be a commutative square in $\cC$ and 
let $(i_a,q_q,u_a)$ be a factorization of $a$. 
Then there exists a factorization $(i_b,q_b,u_b)$ of $b\colon y\to y'$ 
and a morphism $j\colon u_a\to u_b$ which makes the diagram below 
commutative and the canonical morphism $u_a\sqcup_x y\to u_b$ 
induced from the universal 
property of pushout diagram is a cofibration.
\begin{equation}
\label{eq:factorization diagram}
\xymatrix{
x \ar[r]^i \ar@{>->}[d]_{i_a} & y \ar@{>->}[d]^{i_b}\\
u_a \ar[r]^j \ar[d]_{q_a} & u_b \ar[d]^{q_b}\\
y \ar[r]_{i'} & y'.
}
\end{equation}
Moreover assume that $i$ is a cofibration. 
Then we can take $j$ as a cofibration.
\end{lem}

\begin{proof}
Let $(i'',q_b,u_b)$ be a factorization of $u_a\sqcup_x y\to y'$ 
the induced morphism from the universal property of pushout. 
We write $i_b$ and $j$ for the compositions 
$y\rinf u_a\sqcup_x y \overset{i''}{\rinf} u_b$ and 
$u_a  \to u_a\sqcup_x y \overset{i''}{\rinf} u_b$ respectively. 
Moreover if $i$ is a cofibration, 
then the induced morphism $u_a\rinf u_a\sqcup_xy$ is a cofibration. 
Therefore $j$ is a cofibration. 
\end{proof}

\begin{lem}
\label{lem:pushout of factorization}
Let 
\begin{equation}
\label{eq:pushout ladder}
\xymatrix{
y \ar[d]_b & x \ar@{>->}[l]_{i}  \ar[r]^p \ar[d]_a & z \ar[d]^c\\
y' & x \ar@{>->}[l]^{i'}  \ar[r]_{p'} & z'
}
\end{equation}
be a commutative diagram in $\cC$ with $i$ and $i'$ cofibrations and 
let $(i_a,q_a,u_a)$, $(i_b,q_b,u_b) $ and 
$(i_c,q_c,u_c)$ be factorizations of $a$, $b$ and $c$ 
respectively 
and let $i''\colon u_a\to u_b$ and $p''\colon u_a\to u_c$ 
be morphisms which makes diagram below commutative.
\begin{equation}
\label{eq:pushout factorization}
\xymatrix{
y \ar[d]_{i_b} & x \ar@{>->}[l]_{i}  \ar[r]^p \ar[d]_{i_a} & z \ar[d]^{i_c}\\
u_b \ar[d]_{q_b} & u_a \ar@{>->}[l]_{i''}  \ar[r]^{p''} \ar[d]_{q_a} & 
u_c \ar[d]^{q_c}\\
y' & x \ar@{>->}[l]^{i'}  \ar[r]_{p'} & z'.
}
\end{equation}
Assume that $y\sqcup_xu_a \to u_b$ 
the induced morphism from the universal property of pushout 
is a cofibration. 
Then the triple $(i_b\sqcup_{i_a}i_c,q_b\sqcup_{q_a}q_c,u_b\sqcup_{u_a}u_c)$ 
is a factorization of $b\sqcup_a c$.
\end{lem}

\begin{proof}
$i_b\sqcup_{i_a}i_c$ is a cofibration by \cite[1.1.2]{Wal85}.
$q_b\sqcup_{q_a}q_c$ is a morphism in $w_{\bB}$ by the gluing axiom. 
$b\sqcup_a c=(q_b\sqcup_{q_a}q_c)(i_b\sqcup_{i_a}i_c)$. 
Thus the triple $(i_b\sqcup_{i_a}i_c,q_b\sqcup_{q_a}q_c,u_b\sqcup_{u_a}u_c)$ is a factorization of $b\sqcup_a c$.
\end{proof}

\begin{df}[\bf $\cS$-weak equivalences]
\label{df:S-weak equivalence}
Let $\cS$ be a $w$-left thick subcategory of $\cC$. 
Then for a morphism $x\onto{f}y$ in $\cC$, 
the following two conditions are equivalent by 
Lemma~\ref{lem:well-def of S-weak equiv} below and 
in this case, we say that 
$f$ is an {\it $(\cS,w)$-weak equivalence} or simply 
{\it $\cS$-weak equivalence}.
\begin{itemize}
\item
There exists a factorization $(i_f,p_f,u_f)$ of $f$ 
such that $u_f/x$ is in $\cS$. 

\item
For any factorization $(i_f,p_f,u_f)$ of $f$, 
$u_f/x$ is in $\cS$. 
\end{itemize}
We denote the class of all $(\cS,w)$-weak equivalences in $\cC$ 
by $w_{\cS,\cC,w}$ or shortly $w_{\cS}$.
\end{df}

\begin{lem}
\label{lem:well-def of S-weak equiv}
Let $\cS$ be a $w$-left thick subcategory of $\cC$. Then
\begin{enumerate}
\enumidef
\item
Let $x\overset{i}{\rinf}y\overset{j}{\rinf}z$ be a pair of 
composable cofibrations in $\cC$ such that $z/y$ is in $\cS$. 
Then $y/x$ is in $\cS$ if and only if $z/x$ is in $\cS$. 
\item
If $x\overset{i}{\rinf}y$ is a cofibration and a morphism in $w$, 
then $y/x$ is in $\cC^w$.
\item
In the commutative diagram below, assume that 
both $i$ and $i'$ are cofibrations and both $p$ and $p'$ are 
morphisms in $w$,
$$
{\footnotesize{
\xymatrix{
x \ar@{>->}[r]^i \ar@{>->}[d]_{i'} & z \ar[d]^p\\
z' \ar[r]_{p'} & y,
}}}
$$
then 
$z/x$ is in $\cS$ if and only if $z'/x$ is in $\cS$.
\end{enumerate}
\end{lem}

\begin{proof}
$\mathrm{(i)}$ 
By considering a cofibrations sequence 
$y/x\rinf z/x \rdef z/y$, we obtain the result.  

\sn
$\mathrm{(ii)}$ 
By the commutative diagram below, it turns out that 
the canonical morphism $y/x \to 0$ is in $w$ by the gluing axiom. 
$$
{\footnotesize{
\xymatrix{
x \ar@{>->}[r]^i \ar@{>->}[d]_i & y \ar@{=}[d] \ar@{->>}[r] & y/x \ar[d]\\
y \ar@{=}[r] & y \ar[r] & 0.
}}}
$$
Then since $w$ satisfies two out of three property for compositions, 
it turns out that $y/x$ is in $\cC^w\subset \cS$. 

\sn
$\mathrm{(iii)}$
By the factorization axiom, 
there exists a factorization 
$(i_{p\sqcup p'},q_{p\sqcup p'},u)$ of $p\sqcup p'\colon z\sqcup z'\to y$. 
Then by the two out of three property for $w$, 
the composition $z\to z\sqcup z' \onto{i_{p\sqcup p'}} u$ 
is a cofibration and a morphism in $w$. 
Thus by $\mathrm{(ii)}$, $u/z$ is in $\cS$. 
Then by $\mathrm{(i)}$, $z/x$ is in $\cS$ if and only if $u/x$ is in $\cS$. 
By symmetry, it is also equivalent to the condition that $z'/x$ is in $\cS$.  
\end{proof}

\begin{prop}
\label{prop:adm and Serre}
Let $\cS$ be a $w$-thick subcategory of $\cC$ and let 
$u$ be an admissible class of morphisms in $\cC$ which contains $w$. 
Then 
\begin{enumerate}
\enumidef
\item
$w_{\cS}$ is an admissible class of morphisms in $\cC$ which contains 
$w$. 
\item
We have the equalities
\begin{equation}
\label{eq:CwS=S}
\cC^{w_{\cS}}=\cS,
\end{equation}
\begin{equation}
\label{eq:wCu=u}
w_{\cC^u}=u.
\end{equation}
\end{enumerate}
\end{prop}

\begin{proof}[Proof of Proposition~\ref{prop:adm and Serre}]
$\mathrm{(i)}$ 
For a morphism $f\colon x\to y$ in $w$, 
the triple $(\id_x,f,x)$ is a factorization of $f$ such that $x/x\isoto 0$ is in $\cS$. 
Thus $f$ is in $w_{\cS}$. 
In particular, $w_{\cS}$ contains all isomorphisms. 

We consider 
the commutative diagram of 
cofibration sequences $\mathrm{(\ref{eq:cofib seq})}$. 
Let $(i_a,q_a,u_a)$ be a factorization of $a$. 
Then by Lemma~\ref{lem:func factorization} and 
Lemma~\ref{lem:pushout of factorization}, 
there exists factorizations $(i_b,q_b,u_b)$ and $(i_c,q_c,u_c)$ of 
$b$ and $c$ respectively 
and a cofibration sequence $u_a\rinf u_b \rdef u_c$ such that the induced morphisms $u_a\sqcup_xy \to u_b$ 
is a cofibration and 
they make the diagram below commutative: 
$$
\footnotesize{
\xymatrix{
x \ar@{>->}[r] \ar@{>->}[d]_{i_a} & y \ar@{->>}[r] \ar@{>->}[d]_{i_b} & 
y/x \ar@{>->}[d]^{i_c}\\
u_a \ar@{>->}[r] \ar[d]_{q_a} & u_b \ar@{->>}[r] \ar[d]_{q_b} & 
u_c \ar[d]^{q_c}\\
x' \ar@{>->}[r] & y' \ar@{->>}[r] & (y/x).
}}
$$
Then the sequence $u_a/x\rinf u_b/y \rdef u_c/(y/x)$ is 
a cofibration sequence by \cite[1.1.2]{Wal85}. 
Hence if two of $u_a/x$, $u_b/y$ and $u_c/(y/x)$ are in $\cS$, 
then the third one is also in $\cS$. Namely 
if two of $a$, $b$ and $c$ are in $w_{\cS}$, then 
the third one is also in $w_{\cS}$.

Next let $x\onto{a}y\onto{b}z$ be composable morphisms in $\cC$. 
Let $(i_a,q_a,u_a)$ be a factorization of $a$. 
Let $(i',q_{ba},u_{ba})$ be a factorization of $bq_a\colon u_a\to z$. 
Then we set $i_{ba}:=i'i_a$ and it turns out that 
the triple $(i_{ba},q_{ba},u_{ba})$ 
is a factorization of $ba\colon x\to z$. 
By Lemma~\ref{lem:func factorization}, 
there exists a factorization 
$(i_{ba\sqcup\id_y},q_{ba\sqcup\id_y},u_{ba\sqcup\id_y})$ of 
$ba\sqcup\id_y\colon x\sqcup y\to z\sqcup y$ and 
a cofibration $j\colon u_{ba}\rinf u_{ba\sqcup\id_y}$ such that 
$u_{ba}\sqcup_x(x\sqcup y)\rinf u_{ba\sqcup\id_y}$ 
the induced morphism from the universal property 
of pushout is a cofibration. 
Then by Lemma~\ref{lem:pushout of factorization}, there exists 
factorizations $(i_b,q_b,u_b)$, $(i_{\id_y},q_{\id_y},u_{\id_y})$ and 
$(i_{ba\sqcup\id_y},q'_{ba\sqcup\id_y},u_{ba\sqcup\id_y})$ of 
$b$, $\id_y$ and $ba\sqcup\id_y$ respectively 
and cofibration sequences 
$u_a\overset{ji'}{\rinf}u_{ba\sqcup\id_y}\rdef u_b$ and 
$u_{ba}\overset{j}{\rinf}u_{ba\sqcup\id_y}\rdef u_{\id_y}$ 
such that 
they makes the diagrams below commutative:
$$
\xymatrix{
x \ar@{>->}[r] \ar@{>->}[d] & x\sqcup y \ar@{>->}[d] \ar@{->>}[r] & 
y \ar@{>->}[d]\\
u_a \ar@{>->}[r] \ar[d] & u_{ba\sqcup\id_y} \ar@{->>}[r] \ar[d] & 
u_{b} \ar[d]\\
y \ar@{>->}[r] & z\sqcup y \ar@{->>}[r] & z
}\ \ 
\xymatrix{
x \ar@{>->}[r] \ar@{>->}[d] & x\sqcup y \ar@{>->}[d] \ar@{->>}[r] & 
y \ar@{>->}[d]\\
u_{ba} \ar@{>->}[r] \ar[d] & u_{ba\sqcup\id_y} \ar@{->>}[r] \ar[d] & 
u_{\id_y} \ar[d]\\
z \ar@{>->}[r] & z\sqcup y \ar@{->>}[r] & y.
}
$$
By \cite[1.1.2]{Wal85}, the following sequences are cofibration 
sequences
\begin{equation}
\label{eq:u cof seq 1}
u_a/x \rinf u_{ba\sqcup\id_y}/(x\sqcup y) \rdef u_b/y,
\end{equation}
\begin{equation}
\label{eq:u cof seq 2}
u_{ba}/x \rinf u_{ba\sqcup\id_y}/(x\sqcup y) \rdef u_{\id_y}/y.
\end{equation}
Since $\id_y$ is a cofibration and a morphism in $w$, 
$u_b/y$ is in $\cC^w\subset \cS$ by 
Lemma~\ref{lem:well-def of S-weak equiv} $\mathrm{(iii)}$. 
Thus in the cofibration sequence $\mathrm{(\ref{eq:u cof seq 2})}$, 
$u_{ba}/x$ is in $\cS$ if and only if $u_{ba\sqcup\id_y}/(x\sqcup y)$ 
is in $\cS$. 
Hence if two of $u_a/x$, $u_{ba}/x$ and $u_b/y$ are in $\cS$, 
then the third one is also in $\cS$. 
Namely $w_{\cS}$ satisfies the two out of three property for compositions.

\sn
$\mathrm{(ii)}$ 
We show the equality $\mathrm{(\ref{eq:CwS=S})}$. 
Let $x$ be an object in $\cC$ and 
let $(i_0,q_0,u_0)$ be a factorization of the canonical morphism $x\to 0$. 
Then since $p_0\colon u_0 \to 0$ is in $w$, $u_0$ is in $\cC^w\subset \cS$. 
Thus in the cofibration sequence $x\overset{i_0}{\rinf} u_0\rdef u_0/x$, 
$x$ is in $\cS$ if and only if $u_0/x$ is in $\cS$. 
Thus $x$ is in $\cS$ if and only if $x\to 0$ is in $w_{\cS}$ 
if and only if $x$ is in $\cC^{w_{\cS}}$.

Next we prove the equality $\mathrm{(\ref{eq:wCu=u})}$. 
Let $x\onto{f}y$ be a morphism in $\cC$ and let 
$(i_f,p_f,u_f)$ be a factorization of $f$. 
Then $f$ is in $w_{\cC^u}$ if and only if $u_f/x$ is in $\cC^u$. 
In the commutative diagram below, 
$$
{\footnotesize{
\xymatrix{
x \ar@{>->}[r]^{i_f} \ar[d]_f & u_f \ar@{->>}[r] \ar[d]_{p_f} & u_f/x \ar[d]\\
y \ar@{=}[r] & y \ar[r] & 0,
}}}
$$
$f$ is in $u$ if and only if 
the canonical morphism $u_f/x\to 0$ is in $u$ and the last condition 
is equivalent to the condition that $0\to u_f/x$ is in $u$ by 
the two out of three property of $u$.  
\end{proof}

\begin{cor}
\label{cor:w-thick is w-closed}
Let $\cS$ be a $w$-thick subcategory of $\cC$. 
Then $\cS$ is $w$-closed.
\end{cor}

\begin{proof}
Since $w_{\cS}$ contains $w$, $\cC^{w_{\cS}}$ is $w$-closed 
by Lemma~\ref{lem:ex of w-Serre subcategory}. 
By the equality $\mathrm{(\ref{eq:CwS=S})}$, 
$\cS$ is also $w$-closed. 
\end{proof}

\section{Puppe fibration sequences for connective $K$-theory}

Let $f\colon\bA=(\cA,w_{\bA})\to \bB=(\cB,w_{\bB})$ be a morphism of 
Waldhausen categories. 
We write $\im f$ for the class of those objects $x$ in $\cB$ such that 
there exists an object $y$ in $\cA$ and an isomorphism $x\isoto f(y)$. 

\begin{df}[\bf Admissible morphisms]
\label{df:Admissible morphism}
\begin{enumerate}
\enumidef
\item
We say that $f$ is {\it preadmissible} if 
$\im f$ is a $w_{\bB}$-left thick subcategory 
(see Definition~\ref{df:thick subcategory}) 
of $\cB$. 

\item
We say that $f$ is {\it semi-admissible} if it is 
preadmissible and for any cofibration $i\colon f(x)\rinf f(y)$ in $\cB$, 
there exists a cofibration $i'\colon x\rinf y$ in $\cA$ such that 
$f(i')=i$.

\item
We say that $f$ is {\it admissible} if it satisfies the following four conditions hold.
\begin{enumerate}
\enumiidef
\item
$f$ is preadmissible.
\item
The association
$f\colon \Ob\cA \to \Ob\cB$ on the classes of objects are injective.
\item
For any cofibration sequence 
$f(x)\overset{i}{\rinf} f(y) \overset{p}{\rdef} f(z)$, 
there exists a unique cofibration sequence 
$x\overset{i'}{\rinf} y \overset{p'}{\rdef}z$ 
such that $f(i')=i$ and $f(p')=p$. 
\end{enumerate}
\end{enumerate}
\end{df}

\begin{lem}
\label{lem:Im f closed cobase change}
Let $f\colon \bA\to \bB$ be a semi-admissible morphism of 
Waldhausen categories. 
Then $\im f$ is closed under co-base change with respect to cofibrations. 
Namely let $z \overset{j}{\linf}x \to y$ be morphisms 
with $j$ cofibrations and $x$, $y$ and $z\in \im f$. 
Then $z\sqcup_x y$ is also in $\im f$. 
In particular $\im f$ is a $w_{\bB}$-thick subcategory of $\cB$. 
\end{lem}

\begin{proof}
Without loss of generality, 
we shall assume that there exists morphisms 
$x'\overset{j'}{\rinf} z'$ whose image of $f$ is $x\overset{j}{\rinf}z$. 
Then in the commutative diagram of cofibration sequences
$$
\xymatrix{
f(x') \ar@{>->}[r]^{j} \ar[d] & f(z') \ar@{->>}[r] \ar[d] & 
f(z'/x') \ar@{=}[d]\\
z \ar@{>->}[r] & z\sqcup_x y \ar@{->>}[r] & f(z'/x'),
}
$$
$z$ and $f(z'/x')$ are in $\im f$. 
Thus $z\sqcup_x y$ is also in $\im f$. 
\end{proof}

The rest of this section, 
we assume that $\bB$ is strictly derivable 
(see Definition~\ref{df:strictly derivable}) and $f$ is preadmissible. 

\begin{df}
\label{df:f-weak equivalences}
We write $w_f$ for $w_{\im f}$ the class of 
$\im f$-weak equivalences (see Definition~\ref{df:S-weak equivalence}). 
We call morphisms in $w_f$ {\it $f$-weak equivalences}. 
\end{df}

\begin{lem}
\label{lem:w_f property}
Let $f\colon \bA\to \bB$ be 
a semiadmissible morphism of Waldhausen categories. 
Then
\begin{enumerate}
\enumidef
\item
$w_f$ is an admissible class of $\cB$ which contains $w_{\bB}$.

\item
$w_f$ satisfies the gluing axiom. 
\end{enumerate}
\end{lem}

\begin{proof}
$\mathrm{(i)}$ follows from Proposition~\ref{prop:adm and Serre}.

\sn
$\mathrm{(ii)}$ 
We consider the commutative diagram 
$\mathrm{(\ref{eq:pushout ladder})}$ in $\cB$ with $i$ and $i'$ 
cofibrations and $a$, $b$ and $c$ $w_f$-weak equivalences. 
There exists $(i_a,q_a,u_a)$, $(i_b,q_b,u_b)$ and $(i_c,q_c,u_c)$ 
factorizations of $a$, $b$ and $c$ 
respectively and a cofibration $i'\colon u_a\rinf u_b$ and a morphism 
$p''\colon u_a\to u_c$ 
which makes the diagram $\mathrm{(\ref{eq:pushout factorization})}$ 
commutative and $y\sqcup_xu_a\to u_b$ 
the induced morphism is a cofibration by 
Lemma~\ref{eq:commutative square}. 
Then by Lemma~\ref{lem:pushout of factorization}, 
the triple 
$(i_b\sqcup_{i_a}i_c,q_b\sqcup_{q_a}q_c,u_b\sqcup_{u_a}u_c)$ is a factorization of $b\sqcup_a c$ and 
the commutative square 
$$\xymatrix{
u_a/x \ar@{>->}[r] \ar[d] & u_b/y \ar[d]\\
u_c/z \ar@{>->}[r] & (u_c\sqcup_{u_a}u_b)/(z\sqcup_x y)
}$$
is a pushout diagram by \cite[1.1.2.]{Wal85}. 
Since $u_a/x$, $u_b/y$ and $u_c/z$ are in $\im f$, 
$(u_c\sqcup_{u_a}u_b)/(z\sqcup_x y)$ is also in $\im f$ 
by Lemma~\ref{lem:Im f closed cobase change}. 
Thus $b\sqcup_a c$ is in $w_f$. 
\end{proof}

\begin{df}[\bf Cone of morphisms of Waldhausen categories]
\label{df:Cone of morphisms of Waldhausen categories}
Let $\bB$ be a strictly derivable 
Waldhausen category and let 
$f\colon\bA\to \bB$ be a semi-admissible morphism of Waldhausen categories. 
In this case, 
the pair $(\cB,w_f)$ is a Waldhausen category by Lemma~\ref{lem:w_f property}. 
We set $\Cone f:=(\cB,w_f)$ and we denote 
the morphism of Waldhausen categories 
$\bB\to \Cone f$ induced from $\id_{\cB}\colon \cB \to \cB$ by 
$\kappa_f\colon \bB\to \Cone f$. 
We call $\Cone f$ the {\it cone associated with $f$}.
\end{df}

By Lemma~\ref{lem:w_f property}, we obtain the following.

\begin{cor}
\label{cor:Cone strictly derivable}
In the notation 
Definition~\ref{df:Cone of morphisms of Waldhausen categories}, 
$\Cone f$ is a strictly derivable Waldhausen category.
\qed 
\end{cor}

\begin{thm}[\bf Puppe sequences for connective $K$-theory]
\label{thm:Puppe sequences for connective K-theory}
Let $\bB$ be a strictly derivable Waldhausen category and 
let $f\colon\bA\to \bB$ be 
an admissible morphism of Waldhausen categories. 
Then the sequence $\bA\onto{f}\bB\onto{\kappa_f}\Cone f$ induces a fibration sequence of spectra 
$K(\bA)\onto{K(f)}K(\bB)\onto{K(\kappa_f)}K(\Cone f)$ on $K$-theory.
\end{thm}

\begin{proof}
We define $\cB(-w_{\bB})$ to be a simplicial subcategory of 
$\cB^{[-]}$ by sending a totally ordered set $[m]$ to 
$\cB(m,w_{\bB})$ where $\cB(m,w_{\bB})$ is the full subcategory 
of $\cB^{[m]}$ the functor category from $[m]$ to $\cB$ 
consisting of those functors $x\colon [m]\to \cB$ such that 
for any $0\leq i\leq j\leq m$ in $[m]$, 
$x(i\leq j)$ is in $w_{\bB}$. 
A morphism $f\colon x\to y$ in $\cB(m,w_{\bB})$ is a 
{\it level weak equivalence} (respectively {\it level cofibration}) 
if $f_i\colon x_i\to y_i$ is in $w_{\bB}$ (respectively 
a cofibration) for all $0\leq i\leq m$. 
Then we can make $\cB(-,w_{\bB})$ into a 
simplicial Waldhausen category. 

We denote the class of all cofibrations in $w_{f}$ by 
$\bar{w}_{f}$. 
By \cite[1.5.7.]{Wal85}, there exists 
a fibration sequence
$$w_{\bA}S_{\cdot}\cA \to w_{\bB}S_{\cdot}\cB\to wS_{\cdot}S_{\cdot}f.$$
As in \cite[p.352]{Wal85}, 
the third bisimplicial category $wS_{\cdot}S_{\cdot}f$ 
can be identified with $\bar{w}_fS_{\cdot}\cB(-,w_{\bB})$ by 
admissibility of $w_{\bB}$ and $f$ and as in the original proof of 
fibration sequences in \cite{Wal85}, there exists a commutative diagram 
of (bi)simplicial categories
$$\xymatrix{
w_{\bB}S_{\cdot}\cB \ar@{=}[dd] \ar[r] & \bar{w}_fS_{\cdot}\cB(-,w_{\bB}) 
\ar[d]^{\textbf{I}}\\ 
& w_f S_{\cdot}\cB(-,w_{\bB})\\
w_{\bB}S_{\cdot}\cB \ar[r]_{\kappa_f} & w_fS_{\cdot}\cB \ar[u]_{\textbf{II}}
}
$$
where the maps $\textbf{I}$ and $\textbf{II}$ are homotopy equivalences 
by \cite[A.11.]{Sch06} and \cite[1.6.5.]{Wal85} respectively. 
Hence we obtain the result. 
\end{proof}

\section{Unbounded filtered objects}
\label{sec:unbounded filtered objects}

In this section, 
let the pair $\bC:=(\cC,w)$ be 
a strictly derivable Waldhausen category 
(see Definition~\ref{df:strictly derivable}).

\begin{df}[\bf Filtered objects]
\label{df:filtered objects}
Let $\bbN$ be the linearly ordered set of all natural numbers 
with the usual linear order. 
As usual we regard $\bbN$ as a category and 
we denote the category of functors and natural transformations 
from $\bbN$ to $\cC$ by $\cC^{\bbN}$. 
We say that an object $x$ in $\cC^{\bbN}$ is 
a {\it filtered object} ({\it in $\cC$}) 
if $x(n\leq n+1)$ is a cofibration for all $n\in\bbN$. 
We denote the full subcategory of $\cC^{\bbN}$ consisting of 
filtered objects in $\cC$ by $F\cC$. 
For a filtered object $x$ and for a natural number $n$, 
we write $x_n$ and $i^x_n$ for an object $x(n)$ in $\cC$ 
and a morphism $x(n\leq n+1)$ in $\cC$ respectively.

Let $f\colon x\to y$ be a morphism in $F\cC$. 
We say that $f$ is a {\it level weak equivalence} 
if for each natural number $n$, $f_n$ is in $w$. 
We denote the class of all level weak equivalences by ${lw}_{F\cC}$ 
or simply $lw$. 
We say that $f$ is a {\it Reedy cofibration} 
if for each natural number $n$, $f_n$ and the induced 
morphism $x_n\sqcup_{x_{n-1}}y_{n-1}$ are cofibrations.
We can make $F\cC$ into a category of cofibrations by declaring 
the class of all Reedy cofibrations to be the class of cofibrations in $F\cC$ 
(\cf \cite[1.1.4.]{Wal85}). 
\end{df}

We will define several operations on $F\cC$. 

\begin{df}[\bf Truncation functor]
\label{df:truncation}
For a natural number $n$, 
the inclusion functors $F_{\geq n}\cC\to F\cC$ and $F_{\leq n}\cC\to F\cC$ admit right adjoint functors 
$\sigma_{\geq n}\colon F\cC\to F_{\geq n}\cC$ and 
$\sigma_{\leq n}\colon F\cC\to F_{\leq n}\cC$ respectively. 
Explicitly, for an object $x$ and a morphism $f\colon x\to y$ in $F\cC$, 
we set
$$
{(\sigma_{\leq n}x)}_k:=
\begin{cases}
x_k & \text{if $k\leq n$}\\
x_n & \text{if $k\geq n$}
\end{cases},\ \ 
i^{\sigma_{\leq n}x}_k:=
\begin{cases}
i^x_k & \text{if $k\leq n-1$}\\
\id_{x_n} & \text{if $k\geq n$}
\end{cases}
$$
$$
{(\sigma_{\geq n}x)}_k:=
\begin{cases}
x_k & \text{if $k\geq n$}\\
0 & \text{if $k< n$}
\end{cases},\ \ 
i^{\sigma_{\geq n}x}_k:=
\begin{cases}
i^x_k & \text{if $k\geq n$}\\
0 & \text{if $k\leq n-1$}.
\end{cases}
$$
$$
{(\sigma_{\leq n}f)}_k:=
\begin{cases}
f_k & \text{if $k< n$}\\
f_n & \text{if $k\geq n$}
\end{cases},\ \ 
{(\sigma_{\geq n}f)}_k:=
\begin{cases}
f_n & \text{if $k\geq n$}\\
0 & \text{if $k\leq n-1$}.
\end{cases}
$$
There are adjunction morphisms $\sigma_{\geq n}x\to x$ and 
$\sigma_{\leq n}x\to x$.
\end{df}

\begin{df}[\bf Degree shift]
\label{df:degree shift}
Let $n$ be a natural number and let $f\colon x\to y$ be a morphism 
in $F\cC$. 
We define $x[n]$ and $f[n]\colon x[n]\to y[n]$ to be a filtered object in $\cC$ and a morphism of filtered objects respectively by setting 
${x[n]}_k:=x_{n+k}$ and $i^{x[n]}_k:=i^x_{n+k}$ and ${f[n]}_k:=f_{n+k}$. 
The association $(-)[n]\colon F\cC\to F\cC$, $x\mapsto x[n]$ gives an 
exact functor. 

There exists a natural transformation 
$\theta\colon\id_{F\cC}\to \id_{F\cC}[1]$. 
Namely for a filtered object $x$ in $\cC$, 
$\theta_x\colon x\to x[1]$ is defined by the formula 
${\theta_x}_k:=i^x_k$ for each natural number $k$. 
For a pair of natural number $a< b$, 
We write $\theta^{[a,b]}\colon\id_{F\cC}[a]\to \id_{F\cC}[b]$ 
for the compositions $\theta[b-1]\theta[b-2]\cdots\theta[a+1]\theta[a]$. 
For a full subcategory $\calD$ of $F\cC$ which is closed under the operation 
$(-)[1]\colon \calD\to\calD$, we write $\Theta_{\calD}$ for the family 
$\{\theta_x\colon x\to x[1]\}_{x\in\Ob\calD}$ of morphisms in $\calD$ indexed 
by the class of objects in $\calD$. 
\end{df}

\begin{lem}
\label{lem:fundamental propety of lw}
\begin{enumerate}
\enumidef
\item
$lw$ is an admissible class of morphisms in $F\cC$. 

\item
The pair $(F\cC,lw)$ satisfies the factorization axiom.
\end{enumerate}
\end{lem}

\begin{proof}
A proof of assertion $\mathrm{(i)}$ is straightforward. 
We will give a proof of assertion $\mathrm{(ii)}$. 
Let $n$ be a natural number and assume that 
$\sigma_{\leq n}f\colon\sigma_{\leq n}x\to \sigma_{\leq n}y$ 
admits a factorization 
$(i_{\sigma_{\leq n}f},p_{\sigma_{\leq n}f},u_{\sigma_{\leq n}f})$. 
Then let a triple $(i',p_{f_{n+1}},u_{f_{n+1}})$ be 
a factorization of a morphism 
$f_{n+1}\sqcup_{x_n}i_n^yp_{f_n}\colon x_{n+1}\sqcup_{x_n}u_{f_n}\to y_{n+1}$ 
and 
we denote the compositions $x_{n+1}\rinf x_{n+1}\sqcup_{x_n}u_{f_n}\overset{i'}{\rinf}u_{f_{n+1}}$ and $u_{f_n}\to x_{n+1}\sqcup_{x_n}u_{f_n}\overset{i'}{\rinf}u_{f_{n+1}}$ 
by $i_{f_{n+1}}$ and $i_n^{u_f}$ respectively. 
Then the pair of triples $(i_{\sigma_{\leq n}f},p_{\sigma_{\leq n}f},u_{\sigma_{\leq n}f}) $ and $(i_{f_{n+1}},p_{f_{n+1}},u_{f_{n+1}})$ 
give a factorization of 
$\sigma_{\leq n+1}f\colon \sigma_{\leq n+1}x\to \sigma_{\leq n+1}$. 
By proceeding induction on $n$, we finally obtain a factorization of $f$.
\end{proof}

\begin{df}[\bf Bounded and weakly bounded filtered objects]
\label{df:bounded and weakly bounded filtered objects}
Let $a\leq b$ be a pair of natural numbers 
and let $x$ be a filtered object in $\cC$. 
We say that {\it $x$ has amplitude contained in $[a,b]$} if 
for any $0\leq k<a$, $x_k=0$ and for any $b\leq k$, 
$x_k=x_b$ and $i^x_k=\id_{x_b}$. 
In this case we write $x_{\infty}$ for $x_b=x_{b+1}=\cdots$. 
Similarly for any morphism $f\colon x\to y$ in $F\cC$ 
between objects which have amplitude contained in $[a,b]$, 
we denote $f_b=f_{b+1}=\cdots$ by $f_{\infty}$. 
We denote the full subcategory of $F\cC$ consisting of those objects 
having amplitude contained in $[a,b]$ by $F_{[a,b]}\cC$. 
We also set $\displaystyle{F_b\cC:=\underset{\substack{a<b\\ (a,b)\in\bbN^2}}{\bigcup}F_{[a,b]}\cC}$ and 
$\displaystyle{F_{b,\geq a}\cC:=\underset{\substack{a<b\\ b\in\bbN}}{\bigcup}F_{[a,b]}\cC}$ 
and call an object in 
$F_b\cC$ a {\it bounded filtered object} (in $\cC$). 
For a natural number $n$, 
we write $F_{\leq n}\cC$ and $F_{\geq n}\cC$ for 
the category 
$F_{[0,n]}\cC$ and 
the full subcategory of $F\cC$ consisting of those objects $x$ 
such that $x_k=0$ for all $0\leq k<n$ respectively. 

Let $f\colon x\to y$ be a morphism in $\cC$. 
We write $\fj(x)$ and $\fj(f)\colon \fj(x)\to \fj(y)$ for the 
object and the morphism in $F_{[0,0]}\cC$ such that 
${\fj(x)}_{\infty}=x$ and ${\fj(f)}_{\infty}=f$. 
We define ${(-)}_{\infty}\colon F_b\cC\to \cC$ and $\fj\colon \cC\to F_b\cC$ to be functors by sending an object $x$ in $F_b\cC$ to 
$x_{\infty}$ and sending an object $x$ in $\cC$ to $\fj(x)$. 
We denote $F_{[0,0]}\cC$ by $\fj(\cC)$ and sometimes 
identify it with $\cC$ via the functor $\fj$. 

We say that a filtered object $x$ in $\cC$ is {\it weakly bounded} 
if there exists a natural number $N$ such that $i^x_n$ is in $w$ 
for all $n\geq N$. 
We write $F_{wb}\cC$ for the full subcategory of 
$F\cC$ consisting of all weakly bounded filtered objects. 
For a full subcategory $\calD$ of $F\cC$, we denote the class of 
level weak equivalences in $\calD$ by $lw|_{\calD}$ or ${lw}_{\calD}$ or 
simply $lw$.
\end{df}

\begin{df}[\bf Stable weak equivalences]
\label{df:stable weak equivalences}
A morphism $f\colon x\to y$ in $F\cC$ is {\it stably weak equivalence} 
if there exists a natural number $N$ such that 
$f_n$ is in $w$ for all $n\geq N$. 
For a full subcategory $\calD$ of $F\cC$, 
we denote the class of all stably weak equivalences in $\calD$ 
by $w_{\st}|_{\calD}$ or simply $w_{\st}$. 
\end{df}

\begin{lem}
\label{lem:fundamental property of wst}
\begin{enumerate}
\enumidef
\item
$w_{\st}$ is an admissible class of morphisms in $F\cC$. 

\item
The exact functor 
$\fj\colon \bC=(\cC,w) \to (F_b\cC,w_{\st})$ induces 
a homotopy equivalence $K(\bC)\to K(F_b\cC;w_{\st})$ of spectra on $K$-theory.
\end{enumerate}
\end{lem}

\begin{proof}
A proof of assertion $\mathrm{(i)}$ is straightforward. 
For $\mathrm{(ii)}$, notice that we have the equality 
${(-)}_{\infty}\cdot\fj=\id_{\cC}$ and there exists a natural weak equivalence 
$\id_{F_b\cC}\to \fj\cdot {(-)}_{\infty}$ with respect to $w_{\st}$. 
Thus $\fj$ induces a homotopy equivalence of spectra on $K$-theory. 
\end{proof}

\begin{lem}
\label{lem:char of wst and wsus}
\begin{enumerate}
\enumidef
\item
$w_{\st}|_{F_{wb}\cC}$ is the smallest admissible class of morphisms in $F_{wb}\cC$ which contains $lw$ and $\Theta_{F_{wb}\cC}$. 

\item
$F_{wb}\cC$ is the smallest $w_{\st}$-prethick subcategory of $F\cC$ which 
contains $\fj(\cC)$.

\item
For a morphism $f\colon x\to y$ in $F_{wb}\cC$ with 
$x\in\Ob F_b\cC$, 
there is a factorization $(i_f,p_f,u_f)$ of $f$ with respect 
to $w_{\st}$ the class of stable weak equivalences such that $u_f$ is 
in $F_b\cC$. 
\end{enumerate}
\end{lem}

\begin{proof}
$\mathrm{(i)}$ 
Let $u$ be a class of morphisms in $F_{wb}\cC$ which contains 
$lw$ and $\Theta_{F_{wb}\cC}$. 
Then we will show that $u$ contains $w_{\st}$. 
First we will show that for any natural number $N$ and for any object $x$ in 
$F_{wb}\cC$, the natural morphism $\sigma_{\geq N}x \to x$ is in $u$. 
Notice that in the left commutative diagram below 
$\theta^{[0,N]}_{\sigma_{\geq N}x}\colon \sigma_{\geq N}x \to \sigma_{\geq N}[N]$, $\theta^{[0,N]}_x\colon x\to x[N]$ and $\id_{x[N]}\colon \sigma_{\leq N}x[N]\to x[N]$ are in $u$. 
Thus 
the natural morphism $\sigma_{\geq N}x\to x$ also in $u$ by the two out of three property. 
$$
{\footnotesize{
\xymatrix{
\sigma_{\geq N}x \ar[r] \ar[d]_{\theta^{[0,N]}_{\sigma_{\geq N}x}} & x \ar[d]^{\theta^{[0,N]}_x}\\
(\sigma_{\geq N}x)[N] \ar[r]_{\id_{x[N]}} & x[N],
}\ \ \ \ 
\xymatrix{
\sigma_{\geq M}x \ar[r]^{\sigma_{\geq M}f} \ar[d] & \sigma_{\geq M}y \ar[d]\\
x \ar[r]_f & y. 
}
}}
$$
Next let $f\colon x\to y$ be a morphism in $w_{\st}$ and let 
$M$ be a natural number such that $f_n\colon x_n\to y_n$ is in $w$ for 
all $n\geq M$. 
Then $\sigma_{\geq M}f\colon \sigma_{\geq M}x\to \sigma_{\geq M}y$ 
is in $lw$ and 
there is the right commutative diagram above 
with $\sigma_{\geq M}x\to x$, $\sigma_{\geq M}y\to y$ and $\sigma_{\geq M}f\colon \sigma_{\geq M}x\to \sigma_{\geq M}y$ are in $u$. 
Thus by the two out of three property of $u$, $f$ is also in $u$.

\sn
$\mathrm{(ii)}$ 
First we will show that $F_{wb}\cC$ is a $w_{\st}$-thick subcategory of 
$F\cC$. 
Let $x\rinf y\rdef z$ be a cofibration sequence in $F\cC$. 
Then for each natural number $n$, 
since $w$ is admissible, if 
two of $i_n^x$, $i_n^y$ and $i_n^{y/x}$ are in $w$, 
then the third one is also in $w$.
$$
\footnotesize{
\xymatrix{
x_n \ar@{>->}[r] \ar@{>->}[d]_{i_n^x} & 
y_n \ar@{->>}[r] \ar@{>->}[d]_{i_n^y} & y_n/x_n 
\ar@{>->}[d]^{i^{y/x}_n}\\
x_{n+1} \ar@{>->}[r] & y_{n+1} \ar@{->>}[r] & y_{n+1}/x_{n+1}. 
}}
$$
Namely if two of $x$, $y$ and $y/x$ are in $F_{wb}\cC$, then 
the third one is also in $F_{wb}\cC$. 

Next let $f\colon x\to y$ be a morphism in $w_{\st}$ and assume that 
$x$ (resp. $y$) is in $F_{wb}\cC$. 
Then there exists a natural number $N$ such that $i^x_n$ (resp. $i^y_n$) 
and $f_n$ are in $w$ for all $n\geq N$. 
Then by considering the commutative diagram below and the two out of 
three property, it turns out that $i_n^y$ (resp. $i_n^x$) 
is also in $w$.
$$
\footnotesize{
\xymatrix{
x_n \ar[r]^{f_n} \ar@{>->}[d]_{i_n^x} & y_n \ar@{>->}[d]^{i_n^y}\\
x_{n+1} \ar[r]_{f_{n+1}} & y_{n+1}.
}
}
$$
Thus $y$ (resp. $x$) is also in $F_{wb}\cC$. 

Next let $\cS$ be a $w_{\st}$-prethick subcategory of $F\cC$ which contains 
$\fj(\cC)$. 
We will show that $\cS$ contains $F_{wb}\cC$. 
For an object $x$ in $F_{wb}\cC$, there is an integer $N$ 
such that $i_n^x$ is in $w$ for all $n\geq N$. 
Then there is a zig-zag sequence $\fj(x_N)\leftarrow \sigma_{\geq N}\fj(x_N) \to x$ of morphisms in $w_{\st}$. 
Since $\cS$ is $w_{\st}$-closed, it contains $x$. 

\sn
$\mathrm{(iii)}$ 
If $y$ is in $F_b\cC$, then by 
the same proof of \ref{lem:fundamental propety of lw} $\mathrm{(ii)}$, 
we can show that 
there exists a factorization 
$(i_f,p_f,u_f)$ 
of $f$ with $u_f\in\Ob F_b\cC$ and $p_f\in lw$. 

In general case, let $N$ be a natural number such that $i_n^x=\id_{x_N}$ 
and $i_n^y\in w$ for any $n\geq N$. 
Then there exists a factorization $x\onto{f'}\sigma_{\leq N}x\onto{p'}y$ of $f$ with $p'\in w_{\st}$. 
Now by applying the argument in the previous paragraph to $f'$, 
we obtain the desired factorization.
\end{proof}

\begin{cor}
\label{cor:C and FwbC}
The inclusion functor $\bC=(\cC,w)\onto{\fj}(F_b\cC,w_{\st})\to(F_{wb}\cC,w_{\st})$ induces a homotopy equivalence 
$K(\bC)\to K(F_{wb}\cC;w_{\st})$ 
of $K$-theory.
\end{cor}

\begin{proof}
By approximation theorem in \cite[A.2]{Sch06} and \ref{lem:char of wst and wsus} $\mathrm{(iii)}$, 
the inclusion functor $(F_b\cC,w_{\st})\rinc (F_{wb}\cC,w_{\st})$ induces a homotopy equivalence 
of spectra on $K$-theory. 
Thus combining with \ref{lem:fundamental property of wst} 
$\mathrm{(ii)}$, 
we obtain the result.  
\end{proof}

\begin{df}
\label{df:FC and SC}
We denote the smallest admissible class of morphisms 
in $F\cC$ which contains $lw$ and $\Theta_{F\cC}$ by $w_{\sst}$ and call it 
the {\it class of strongly stably weak equivalences in $F\cC$}. 
We write $w_{\sus}$ for 
the admissible class 
$w_{F_{wb}\cC}$ of morphisms 
in $F\cC$ which corresponds to the $w_{\st}$-prethick subcategory 
$F_{wb}\cC$ of $F\cC$ and call it 
the {\it class of suspension weak equivalences in $F\cC$}. 
We write $\cF\bC$ and $\fS\bC$ for the pairs 
$(F\cC,w_{\sst})$ and $(F\cC,w_{\sus})$ and call them 
{\it flasque envelope of $\bC$} and the {\it suspension of $\bC$} 
respectively. 
The naming of $\fF\bC$ and $\fS\bC$ are justified by the following results. 
\end{df}

\begin{prop}
\label{prop:Eilenberg swindle and delooping}
Assume that $\cC$ is essentially small. Then
\begin{enumerate}
\enumidef
\item
{\bf (Eilenberg swindle).}\ \ 
$K(\fF\bC)$ is trivial. 

\item
{\bf (Suspension).}\ \ 
There exists a natural homotopy equivalence of spectra 
$
K(\fS\bC)\isoto \Sigma K(\bC)
$.
\end{enumerate}
\end{prop}

\begin{proof}
$\mathrm{(i)}$ 
We denote the functor $F\cC\to F\cC$ which sends an object $x$ in $F\cC$ to 
$\bigsqcup_{n\geq 0}x[n]$ by $F$. 
Then there exists an equality $F=F[1]\sqcup\id_{F\cC}$ and there 
exists a natural weak equivalence 
$\theta_F\colon F\to F[1]$ with respect to 
$w_{\sst}$. 
Thus $K(F)=K(F[1])$ on $K(\fF\bC)$ and the identity morphism of 
$K(\fF\bC)$ is trivial. 
Hence $K(\fF\bC)$ is trivial.

\sn
$\mathrm{(ii)}$ 
By fibration theorem in \cite[A.3]{Sch06} and \ref{cor:C and FwbC}, 
there exists a fibration sequence 
$K(\bC)\to K(\fF\bC)\to K(\fS\bC)$ 
of spectra. 
Combining with the result $\mathrm{(i)}$, 
we obtain the homotopy equivalence of spectra $K(\fS\bC)\isoto \Sigma K(\bC)$. 
\end{proof}

\begin{df}[\bf Non-connective $K$-theory for strictly derivable Waldhausen categories]
\label{df:nonconnective K-theory}
For an essentially small strictly derivable Waldhausen category 
$\bC=(\cC,w)$, 
by Proposition~\ref{prop:Eilenberg swindle and delooping} $\mathrm{(ii)}$, 
there exists a canonical morphism
\begin{equation}
\label{eq:structure map}
K(\bC)\to \Omega K(\fS\bC).
\end{equation} 
We denote the spectrum $\bbK(\bC)$ to be the sequence of 
spaces $K(\bC)$, $K(\fS\bC)$, $K(\fS^2\bC)$ 
with the structure morphisms given by 
$\mathrm{(\ref{eq:structure map})}$. 
\end{df}

\mn
SATOSHI MOCHIZUKI\\
{\it{DEPARTMENT OF MATHEMATICS,
CHUO UNIVERSITY,
BUNKYO-KU, TOKYO, JAPAN.}}\\
e-mail: {\tt{mochi@gug.math.chuo-u.ac.jp}}\\

\end{document}